\documentclass[%
 aip,
 amsmath,amssymb,
 reprint,%
]{revtex4-1}

\usepackage[english]{babel}
\usepackage{graphicx}% Include figure files
\usepackage{dcolumn}% Align table columns on decimal point
\usepackage{bm}% bold math
\usepackage{mathtools}% for \mathtoolsset{showonlyrefs}, which numbers only used equations
\usepackage[colorlinks=true,linkcolor=blue,urlcolor=blue,citecolor=blue,breaklinks=true]{hyperref}

\usepackage{listings} % for listing code
\usepackage{xcolor}

\definecolor{codegreen}{rgb}{0,0.6,0}
\definecolor{codegray}{rgb}{0.5,0.5,0.5}
\definecolor{codepurple}{rgb}{0.58,0,0.82}
\definecolor{backcolour}{rgb}{0.95,0.95,0.98}

\lstdefinestyle{mystyle}{
    backgroundcolor=\color{backcolour},   
    commentstyle=\color{codegreen},
    keywordstyle=\color{magenta},
    numberstyle=\tiny\color{codegray},
    stringstyle=\color{codepurple},
    basicstyle=\ttfamily\footnotesize,
    breakatwhitespace=false,         
    breaklines=true,                 
    captionpos=b,                    
    keepspaces=true,                 
    showspaces=false,                
    showstringspaces=false,
    showtabs=false,                  
    tabsize=2
}

\usepackage[utf8]{inputenc}
\usepackage[T1]{fontenc}
\usepackage{mathptmx}
\usepackage{etoolbox}
\usepackage{float}

\newcommand{\norm}[1]{\| #1 \|}

\newcommand{\dint}{\mathrm{d}}

\newcommand{\comments}[1]{}
\newcommand{\R}{\mathbb{R}}
\newcommand{\A}{\mathcal{A}}
\newcommand{\x}{\mathbf{x}}
\newcommand{\y}{\mathbf{y}}
\newcommand{\B}{\mathcal{B}}
\newcommand{\dist}{\mathrm{dist}}
\mathtoolsset{showonlyrefs} % numbering equations
\makeatletter
\def\@email#1#2{%
 \endgroup
 \patchcmd{\titleblock@produce}
  {\frontmatter@RRAPformat}
  {\frontmatter@RRAPformat{\produce@RRAP{*#1\href{mailto:#2}{#2}}}\frontmatter@RRAPformat}
  {}{}
}%
\makeatother

\begin{document}

\preprint{AIP/123-QED}

\title{Computing Resilience Measures in Dynamical Systems}
%\title{Numerically estimating resilience concepts in dynamical systems}

\author{Andreas Morr}
\email{andreas.morr@tum.de}
\affiliation{Department of Mathematics, School of Computation, Information and Technology, Technical University of Munich, Germany}
\affiliation{Research Domain IV -- Complexity Science, Potsdam Institute for Climate Impact Research, Germany}

\author{Christian Kuehn}%
\affiliation{Department of Mathematics, School of Computation, Information and Technology, Technical University of Munich, Germany}

\author{George Datseris}%
\affiliation{Department of Mathematics and Statistics, University of Exeter, Exeter, United Kingdom}%

\date{\today}

\begin{abstract}
Resilience broadly describes a quality of withstanding perturbations. Measures of system resilience have gathered increasing attention across applied disciplines, yet existing metrics often lack computational accessibility and generalizability. In this work, we review the literature on resilience measures through the lens of dynamical systems theory and numerical methods. In this context, we reformulate pertinent measures into a general form and introduce a resource-efficient algorithm designed for their parallel numerical estimation. By coupling these measures with a global continuation of attractors, we enable their consistent evaluation along system parameter changes. The resulting framework is modular and easily extendable, allowing for the incorporation of new resilience measures as they arise. We demonstrate the framework on a range of illustrative dynamical systems, revealing key differences in how resilience changes across systems. This approach provides a more global perspective compared to traditional linear stability metrics used in local bifurcation analysis, which can overlook inconspicuous but significant shifts in system resilience. This work opens the door to genuinely novel lines of inquiry, such as the development of new early warning signals for critical transitions or the discovery of universal scaling behaviours. All code and computational tools are provided as an open-source contribution to the DynamicalSystems.jl software library.
\end{abstract}

\maketitle

\begin{quotation}
Resilience concepts have inspired systems analysis methodology across scientific disciplines and have shaped the public discourse on hazards and potential harms. \cite{Naderpajouh2023ResilienceReview, CatauVeres2021ResiliencePublicDiscourse} 
Generally, resilience describes an inherent quality of a system to withstand external influences, and can be thought of as a generalization of ``stability''. \cite{Gruemm1976DefinitionsResilience, Hosseini2016ResilienceMeasuresEconomy, Fraccascia2018ResilienceComplexSystems}
Assessing an often undesirable loss of resilience has become the focus of a growing field of research. \cite{Carpenter2011EWSEcosystemExperiment, Boers2021AMOCEWS, Boulton2022AmazonResilience, Smith2022VegetationResilience, Han2024HumanImpactsResilience} This is particularly important in the context of a changing world. \cite{Moore2022ResilienceClimateChange, Rockstrom2023PlanetBoundaries} Ecological systems can lose resilience under increasing stress from climate change, \cite{Pimm1984Resilience} while financial institutions are regularly probed for their resilience against hypothetical crisis scenarios, \cite{Tahir2025FinancialResilience} to name two examples. In each setting, resilience will have a different meaning and will be measured and estimated through different definitions. \cite{Krakovska2024ResilienceDynamicalSystems} Insofar as the system property of resilience can be quantified by such a measure, its decrease with respect to external changes to the system may point out system vulnerabilities. \cite{Lenton2012CSDClimate, Ditlevsen2023AMOCPrediction} Because such measures have often been defined for specific application contexts, their computation has been an individual endeavour as well. Our framework for the estimation of a broad class of resilience measures in general dynamical systems is going to aid the proliferation and comparability of their applications.
\end{quotation}

\section{Introduction}
In this work, we compare several measures of system resilience that are generally applicable to almost any system regardless of its context. The common basis for these measures is dynamical systems theory, which provides a mathematical framework for modelling and analysing system behaviour. We assume for this work that the time evolution of a system state vector $\x_0\in\R^n$ can be modelled by an ordinary differential equation.

The extent to which any such finite-dimensional description can capture the behaviour of an underlying natural system depends on its complexity.\cite{Givon2004ModelReduction, Hummel2023ClimateModelReduction} Spatially extended systems typically follow partial differential equations and are infinite-dimensional. Often, the rate of change cannot be easily expressed solely as a function of the current state but also of time. Exogenous variables such as noise or parameter changes need to be included, making the right-hand side of the differential equation dependent on time and the system non-autonomous.\cite{Gardiner1985StochAna, Morr2025RedNoise} We choose to focus here on systems that can be effectively represented as autonomous dynamical systems.
As long as the system is stationary, even if non-autonomous, many of the methods and results discussed here translate straightforwardly from the autonomous setting.\cite{Hetzler2023StationaryDynamics} An important subclass of stationary systems is that of periodically forced systems.\cite{Guckenheimer1983NonlinearOscillation, Williamson2016EWSPeriodic} For other non-autonomous cases, whether the concepts carry over has to be examined on a case-by-case basis.
The general autonomous differential equation of interest in this work is
\begin{equation}
    \dot\x:=\frac{\dint \x}{\dint t}=f(\x,p)\label{eq: ODE}
\end{equation}
with the only exogenous variable being a parameter $p\in\R$. We always consider the right-hand side $f:\R^n\times\R\rightarrow\R^n$ to depend continuously on $p$. Furthermore, we demand that for fixed $p$, the initial value problem corresponding to Eq.~\eqref{eq: ODE} and any $\x_0\in\R^n$ admits a unique solution $\x(\cdot,\x_0)$. In this setup, we may compare resilience measures of autonomous systems for any parameter value $p$ and assess their response to parameter changes.

The chosen array of resilience measures is sourced from various fields of applied and mathematical research and intentionally covers a diverse set of resilience concepts. They can be broadly assigned to one of two categories: local and nonlocal resilience measures. The terms engineering and ecological resilience, respectively, are often used interchangeably. \cite{Holling1973ResilienceStabilityEcological, Dakos2022ResilienceMeasures} Locality indicates that the measures are interpreted with respect to the immediate neighbourhood of an attracting system state, often a stable steady state. Nonlocal resilience measures, on the other hand, synthesize information about the entire state space. Perhaps the most popular indicator of changing system resilience, so-called critical slowing down, is directly related to the linear stability of the attractor and is therefore a local measure. This observational characteristic manifests when a dynamical system, perturbed by randomness closely around its stable steady state, loses linear stability. The slowing of the perturbed dynamics thus takes place before generic bifurcations and catastrophic shifts, where linear stability is generally lost. For this reason, critical slowing down has been suggested as an early warning signal for such critical transitions \cite{Wiesenfeld1985EWSNonlinearInstab, Dakos2008CSDClimate, Scheffer2009EWSNature, Lenton2012CSDClimate, Kuehn2013CSDVar}, and many natural systems have been investigated through this lens. \cite{Boers2018DOEWS, Boers2021AMOCEWS, Boers2021CSDGreenlandIS} However, this understanding of system resilience fails to include information about the system's evolution that goes beyond the asymptotic linear response to perturbations. We introduce several such measures that can advance the understanding of different types of resilience in applied systems. In some fields, system ``stability'', ``constancy'' or ``persistence'' are the preferred choice of terminology. \cite{Yi2021ResilienceReview, VanMeerbeek2021OverviewResilience, Datseris2023FrameworkGlobalStability} In this work, we choose to homogenize the language to the single term ``resilience''.

Though the definitions of the considered resilience measures in mathematical terms have been established, the effective and comparative numerical estimation of resilience measures has often not yet been achieved. This is because even with complete knowledge of the system, i.e., its evolution equation \eqref{eq: ODE}, its behaviour needs to be numerically simulated. An appropriately large selection of initial conditions must be probed in order to find attractors of the system and to assess local and nonlocal resilience. \cite{Datseris2023FrameworkGlobalStability} We achieve this by employing the general-purpose Julia library DynamicalSystems.jl~\cite{Datseris2018DynamicalSystemsjl} and extending its code base with functionality to compute and track resilience measures across parameter settings in a parallel manner. We thereby provide an accessible toolkit for researchers in the field of applied dynamical modelling.

Section \ref{sec: resMes} introduces and defines resilience measures and broadly categorises them. Their heuristic interpretation is also discussed. Section \ref{sec: numEst} explains the numerical algorithm employed for the estimation of the resilience measures. In Section \ref{sec: results}, we apply the algorithm to three example systems from climate science and ecology. Their trends in advance of catastrophic regime shifts are assessed and compared. Finally, we discuss the results and their implications for future research in the sphere of systems analysis.

\section{Selected resilience measures}\label{sec: resMes}
To capture a diverse range of resilience concepts and associated measures, we reviewed recent literature employing this terminology. In some cases, concrete definitions in terms of dynamical systems theory were already given. In other cases, a heuristic description needed to be translated into such a definition. The review by Krakovska et al. (Ref.~\onlinecite{Krakovska2024ResilienceDynamicalSystems}) serves as a starting point. Therein, many of the investigated resilience measures are rigorously defined and related to the relevant applied literature. Table \ref{tab: measures} shows a complete list of all measures considered in this work and lists references to previous research and applications in each instance.

The common basis for all further analysis is the autonomous dynamical system of Eq.~\eqref{eq: ODE}. Let $\A\subset\R^n$ be an attractor of this system. This means that for all initial conditions $\x_0$ in some open neighbourhood around $\A$, the trajectories converge to $\A$. Note that we do not presume attractors to be maximal in the sense of set inclusion, and alternative definitions of attractors exist.\cite{Milnor1985ConceptAttractor, Gorodetski1996MinimalStrangeAttractors} Physical knowledge of the system is required to make an astute choice of system attractors that ensures the interpretability of the results.\cite{Schoenmakers2021FunctioningResilience} The attractors we consider are either a stable equilibrium point (steady state), or a stable oscillation (i.e., a limit cycle), or a chaotic attractor (see e.g.~Refs.~\onlinecite{Guckenheimer1983NonlinearOscillation, Strogatz2015NonlinearDynamics, Kuznetsov2004AppliedBifurcationTheory, Datseris2022NonlinearDynamicsJulia} for an introduction). In the context of applied modelling, an attractor can represent a mode of functioning of a physical system. \cite{Radosavljevic2023EcologicalModellingReview, Munch2023EcologicalModelling} Often, some mode of functioning is desired, while a switch to another mode is undesired or even catastrophic. In an example to be discussed further below, an oscillatory state of coexisting predator and prey populations is desired for the sake of biodiversity over the alternative attractors associated with extinction. In this sense, resilience should be measured with respect to each attractor, quantifying its stability concerning potential regime changes.

We define a resilience measure $r$ generally as a function mapping the ODE right-hand side $f(\cdot,p)$ and an attractor $\A$ of that system to a real number:
\begin{equation}
    r^q(f(\cdot,p),\A)\in\R.
\end{equation}
We thereby assign each mode of system functioning a level of resilience measured by $r$. Some measures are specified by additional information $q$, which may also depend on $p$. 

We have assumed the dependence of $f$ on $p$ as continuous. Assume further that there is no dynamic bifurcation in a parameter range $(p_1,p_2)$. This means that the number of attractors is constant, and one can track the evolution of each attractor enumerated by $j\in\mathbb{N}$ as a function $\A_j(p)$ of $p\in(p_1,p_2)$. In that sense, it is possible to speak of the parameter dependence of a resilience measure $r^{q(p)}(f(\cdot,p),\A_j(p))$ for each attractor of a given system.

\subsection{Local Resilience Measures}\label{mes: loc}
\begin{figure}[t]
	\centering
	\includegraphics[width=\columnwidth]{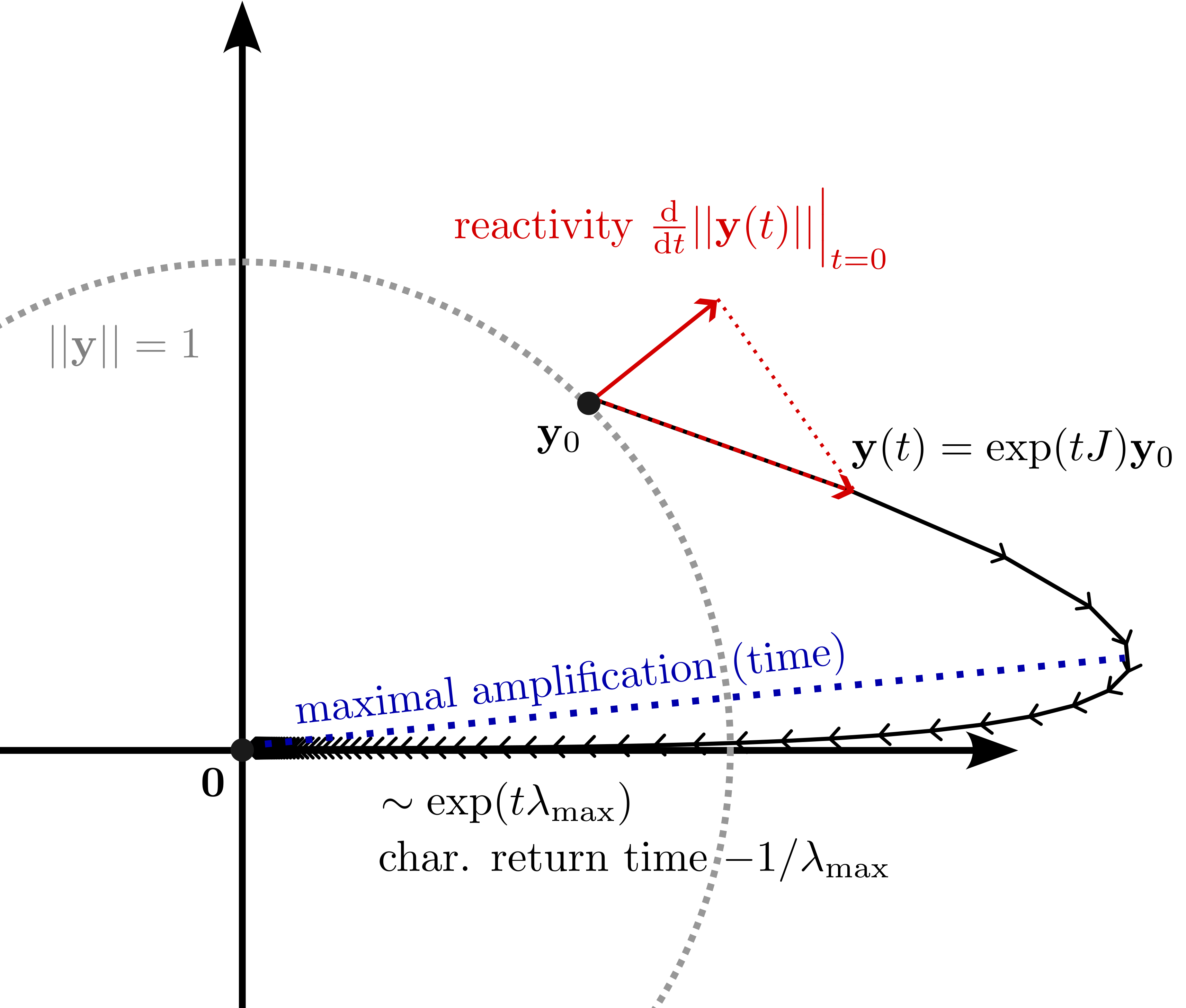}
	\caption{Schematic representation of local resilience measures. The linearised dynamics of a dynamical system around its point attractor are analysed for characteristic return time, reactivity, and maximal amplification (time).}
    \label{fig: localres}
\end{figure}

The class of local resilience measures is concerned with the behaviour of the system in the immediate proximity of the attractor $\A$. If $\A=\{\x^*\}$ is a stable hyperbolic equilibrium point, then the dynamics are characterized by the first-order expansion of $f(\cdot,p)$ around $\x^*$. Therefore, local resilience measures can be seen as functions of the Jacobian matrix $J:=J_\x f(\cdot,p)|_{\x=\x^*}$. The general case of non-point attractors allows for a similar interpretation under the framework of Lyapunov theory or exponential dichotomy, but we will not treat this here. \cite{Dieci2002LyapunovSpectrum, Krakovska2024ResilienceDynamicalSystems} Figure \ref{fig: localres} illustrates the introduced local resilience concepts.

\subsubsection{Characteristic return time}\label{mes: crt}
One generic prerequisite for a system equilibrium point $\x^*$ to be a point attractor is linear stability. Considering the first-order expansion of the dynamics discussed above
\begin{equation}
    \dot{\x}=f(\x,p)\approx J\,(\x-\x^*),
\end{equation}
linear stability is attained if all initial conditions converge to $\x^*$ under the linear equation. For further investigation, we assume that $\x^*=\mathbf{0}$ without loss of generality and define the system 
\begin{equation}
    \dot{\y}=J\,\y
\end{equation}
with solution
\begin{equation}
    \y(t,\y_0)=\exp(tJ)\,\y_0.
\end{equation}
Convergence to $\mathbf{0}$ for all initial conditions $\y_0\in\R^n$ is equivalent to all eigenvalues of $J$ having negative real part. We can bound the rate of convergence in terms of the spectral abcissa $\lambda_\mathrm{max}<0$, the largest real part of any eigenvalue of $J$: There exists some $C>0$ so that for all $\y_0\in \R^n$ and $t\geq0$
\begin{equation}
    \norm{\y(t,\y_0)}\leq C\,\exp(\lambda_\mathrm{max}t)\norm{\y_0},\label{eq: decay}
\end{equation}
where $\norm{\cdot}$ is any norm on $\R^n$, though we will opt for the Euclidean norm $\norm{\cdot}_2$ throughout this work. For this reason, we call 
\begin{equation}
    t_R:=-1/\lambda_\mathrm{max}
\end{equation}
the characteristic return time. Any deviation from the attractor asymptotically shrinks by a factor of at least $1/e$ during each added characteristic return time.

The influence of this concept on resilience research cannot be understated. Linear stability generically vanishes at local bifurcation points. \cite{Kuznetsov2004AppliedBifurcationTheory} This is equivalent to the characteristic return time diverging to infinity. If one is able to assess the linear stability of a natural system from observations, its loss can be interpreted as an approach to a catastrophic shift. \cite{Dakos2008CSDClimate, Scheffer2009EWSNature} As we have alluded to in the Introduction, such indications are often referred to as critical slowing down. \cite{Lenton2012CSDClimate, Dakos2024EWSReview} Popular proxies of linear stability from time series data are variance and lag-1 autocorrelation, \cite{Scheffer2009EWSNature, Ditlevsen2010SignalNoise} and many direct estimators exist. \cite{Weinans2021MultEWSPerformanceComp, Morr2024RedNoiseCSD, Morr2024KramersMoyalEWS} However, numerous theoretical and practical caveats make the application of critical slowing down-based resilience indicators to real-world systems challenging.\cite{Kuehn2013CSDVar, Kuehn2022ColourBlind, Ben-Yami2024TippingTime, Morr2024InternalNoiseInterference, Rietkerk2025AmbiguityEWS}

\subsubsection{Reactivity}\label{mes: reac}
The constant $C$ and the specification that the convergence with rate $\lambda_\mathrm{max}$ in Eq. \eqref{eq: decay} is only asymptotic are important. In fact, the convergence of an initial condition $\y_0$ towards $\mathbf{0}$ does not necessarily commence immediately. There may be a $\y_0$ for which the trajectory first moves away from $\mathbf{0}$ before eventually converging: $\norm{\y(t_1,\y_0)}\geq\norm{\y_0}$. The reactivity 
\begin{equation}
    R_0:=\max_{\norm{\y_0}=1}\left\{\frac{\dint}{\dint t}\norm{\y(t,\y_0)}\Big|_{t=0}\right\}
\end{equation}
represents the largest possible initial growth of trajectories $\y(t,\y_0)$. Restricting to the unit circle is sufficient since the response to any other initial condition is simply scaled by its norm. If the reactivity is smaller than $0$, trajectories immediately move towards $\mathbf{0}$ for any initial condition $\y_0$. If it is larger than $0$, then there exist initial conditions on the unit circle for which the corresponding trajectories first move away from the origin. Reactivity can give an idea of how stabilizing the initial response of the system to a disturbance is, thus capturing a different local resilience concept than the characteristic return time. When it is defined using the Euclidean norm $\norm{\cdot}_2$, it can be directly computed as the spectral abscissa of the matrix $(J+J^\top)/2$. \cite{Arnoldi2016ResilienceMathComparison}

\subsubsection{Maximal amplification (time)}\label{mes: ampl}
To capture the entire finite-time response of the linear system when initialized at any point other than $\mathbf{0}$, we can define the amplification envelope
\begin{equation}\label{eq: amplification}
    \rho(t)=\max_{\norm{\y_0}=1}\left\{\norm{\y(t,\y_0)}\right\}=\norm{\exp(tJ)}_{\mathrm{op}},
\end{equation}
where the right-hand side employs the matrix operator norm.
This definition ensures that the function $\rho$ is larger than the norm of any point on any trajectory started on the unit circle. We define two resilience measures: The maximal amplification
\begin{equation}
    \rho_{\mathrm{max}}:=\max_{t\geq0}\rho(t)
\end{equation}
and the maximal amplification time  
\begin{equation}
    t_{\mathrm{max}}:=\underset{t\geq0}{\mathrm{argmax}}\,\rho(t).
\end{equation}
They capture the worst finite-time response to a disturbance away from $\mathbf{0}$ in terms of the maximal excursion and also specify at which time that excursion reaches its apex. The two quantities can be determined by performing a maximization of the right-hand side of Eq.~\eqref{eq: amplification} over $t\geq 0$.

\subsection{Nonlocal resilience measures}\label{mes: non-loc}
The above local resilience measures are concerned with the system's response to small disturbances away from its attractor. Disturbances need to be small enough for the first-order expansion of the dynamics to possibly be a good approximation. In this subsection, we instead consider the system's response to large disturbances. They may even be large enough for the system to converge to an alternative attractor, a possibility that is disregarded in the local case. In this subsection, we also drop the condition of $\A$ being a point attractor. Instead, $\A$ could be a periodic orbit or a chaotic attractor.

\begin{figure}[t]
	\centering
	\includegraphics[width=\columnwidth]{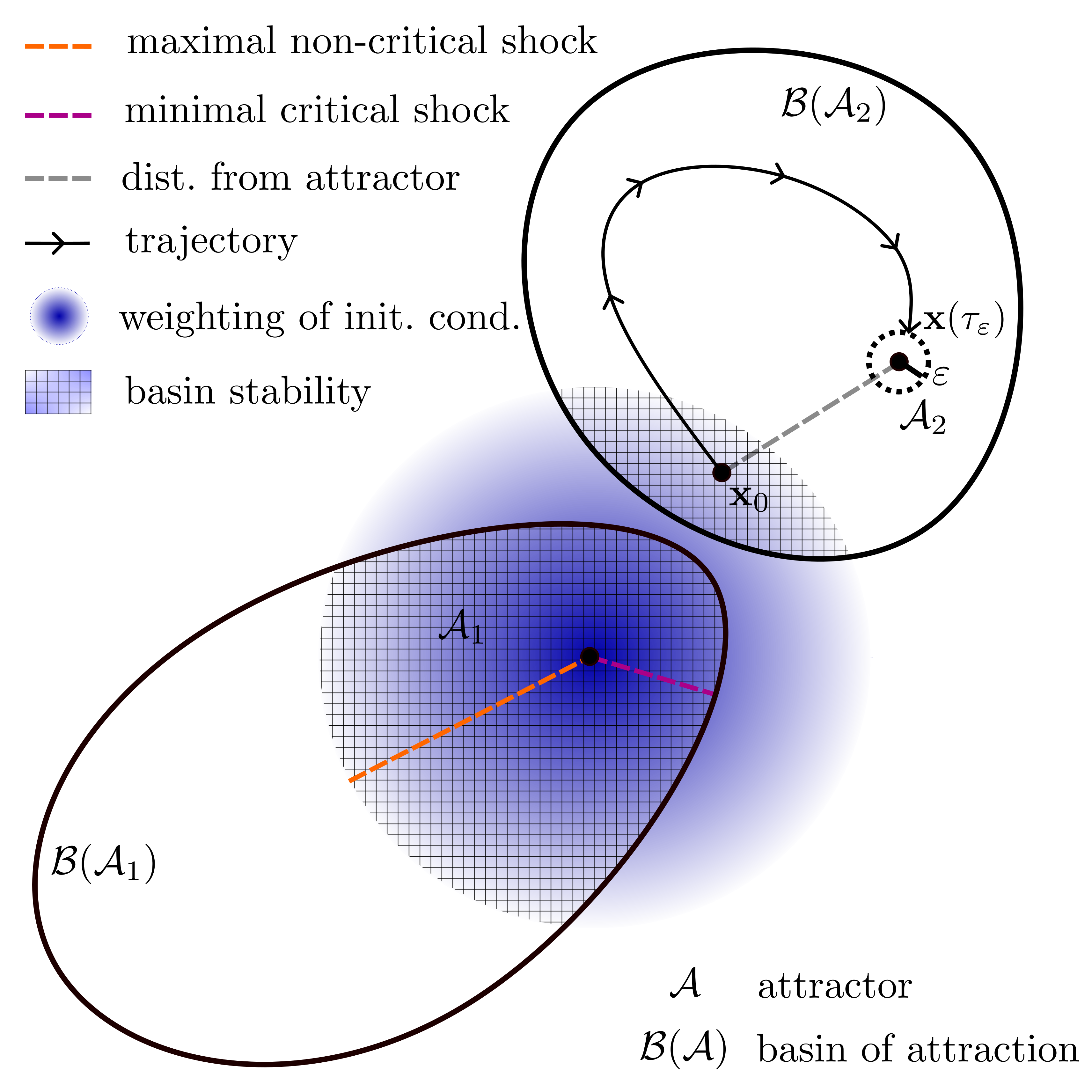}
	\caption{Schematic representation of some nonlocal resilience measures. The size and volume of basins of attraction are gauged with respect to a given probability distribution of initial conditions. This gives the minimal critical shock, maximal critical shock, and basin stability. Taking into account the transient dynamics of trajectories, one obtains convergence time, convergence pace, and finite-time basin stability.}
    \label{fig: nonlocalres}
\end{figure}

Any specific disturbance away from $\A$ may not be of much practical relevance in a given context. Often, the set of possible disturbances clusters around attractors. This could be represented by a probability measure $\nu$ on $\R^n$ encoding the uncertainty of initial conditions $\x_0$. We, therefore, augment the nonlocal resilience measures with a qualifier $q=(\nu)$ to make the results more interpretable in applications. In Figure \ref{fig: nonlocalres}, some of the nonlocal resilience concepts under consideration here are illustrated.

\setcounter{subsubsection}{3}
\subsubsection{Minimal critical shock}\label{mes: mcs}
We call $\B(\A)\subset\R^n$ the basin of attraction associated with $\A$. It is the set of initial conditions $\x_0$ for which $\x(t,\x_0)$ converges to $\A$ as $t\rightarrow\infty$. We are interested in the smallest distance from $\A$ at which we can find an initial condition $\x_0$ not converging to $\A$. It is the minimal magnitude a disturbance must have in order to possibly be critical, i.e., resulting in a new system state asymptotically. However, we only consider those $\x_0$ with positive probability density $w_\nu(\x_0)>0$
\begin{equation}
    s^{(\nu)}_{\mathrm{min}}=\dist(\A,\,(\R^n\setminus\B(\A))|_{w_\nu>0})
\end{equation}
Here, $\dist$ refers to the ordinary distance function of the two sets under a chosen norm $\norm{\cdot}$. More concretely, we use the Euclidean norm $\norm{\cdot}_2$ throughout this work, meaning
\begin{equation*}
    \dist(A,B):=\inf_{a\in A,\,b\in B}\norm{a-b}_2
\end{equation*}
for any $A,\,B\subset\R^n$.
The Euclidean norm may not be suitable for all types of systems. For example, in systems with heterogeneous variables, one may use a weighted Euclidean norm, and in systems with angle variables, one typically uses periodic Euclidean norms. Regardless, if all non-zero probability $\x_0$ converge to $\A$, i.e.~$(\R^n\setminus\B(\A))|_{w_\nu>0}=\varnothing$, then we define the minimal critical shock as infinity.

\subsubsection{Maximal non-critical shock}\label{mes: mncs}
Similarly, we may ask what is the largest magnitude of disturbance from the original attractor $\A$ that could still result in the system converging back to $\A$.
\begin{equation}
    s^{(\nu)}_{\mathrm{max}}=\sup_{\x_0\in\B(\A)|_{w_\nu>0}}\dist(\A,\,\x_0)
\end{equation}

\subsubsection{Basin stability}\label{mes: bs}
The preceding two resilience measures captured the size of $\B(\A)$ in terms of extremal critical or non-critical disturbances. To capture the overall size of $\B(\A)$ with respect to the probable set of disturbances, we may define basin stability to be the set volume under $\nu$:
\begin{equation}
    S^{(\nu)}=\nu(\B(\A)).
\end{equation}
This quantity can be interpreted as the probability of an initial condition $\x_0$ pulled from the prior knowledge distribution $\nu$ to converge to $\A$.

\subsubsection{Convergence time}\label{mes: ct}
In many practical settings, it may be of interest whether a system is within an $\varepsilon$-distance to the attractor. This could signify that the system is close enough to normal functioning. When it is perturbed or initialised at the initial condition $\mathbf{x}_0$, we may be interested in the duration of the transient outside of the $\varepsilon$-distance. We, therefore, define the convergence time of this single initial condition as
\begin{equation}\label{eq: ct}
    \tau^{(\varepsilon)}(\x_0) = \inf\{t\geq0\,|\,\dist(\x(t,\x_0), \A)\leq\varepsilon\}
\end{equation}
To obtain a scalar resilience measure over the entire basin of attraction, one could consider the mean, the maximum, or the median of all initial conditions with respect to $\nu$.

The maximal convergence time will generically only contain non-trivial information if $w_\nu(\x_0)=0$ on the basin boundary $\partial\B(\A)$. That is because, for many smooth systems, the convergence time grows to infinity when approaching the basin boundary, making the maximal convergence time of the entire basin trivially infinite. On the other hand, if an initial condition is closer than $\varepsilon$ to $\A$ to begin with, then $\tau^{(\varepsilon)}(\x_0)$ is trivially zero.

Considering the median convergence time of all non-zero-measure initial conditions has a conceptual and practical advantage. First, it will not be trivially infinite in case the boundary has non-zero measure under the weighting distribution. Second, in a numerical estimation through a finite number of initial conditions, the median will be more robust against outlier values. For the examples discussed in this work, we have chosen to investigate this median quantity rather than the mean or the maximum.

\subsubsection{Convergence pace}\label{mes: cp}
For a fair comparison of the convergence times of different initial conditions $\x_0$, we scale each convergence time by the size of the disturbance to obtain the convergence pace:
\begin{equation}
    \beta^{(\varepsilon)}(\mathbf{x}_0)=\frac{\tau^{(\varepsilon)}(\x_0)}{\dist(\A,\,\x_0)}.
\end{equation}
For $\x_0\in\A$, the convergence pace is set to $0$. This quantity can similarly be synthesized as the mean, maximum, or median over the non-zero probability initial conditions in the respective basin of attraction.

\subsubsection{Finite-time basin stability}\label{mes: ftbs}
If the convergence time $\tau^{(\varepsilon)}(\x_0)$ to get back to an attractor is very large, then the system could also be interpreted to not have recovered from the disturbance $\x_0$ at all in some contexts. For instance, a prolonged transient change of system functioning may already be catastrophic even when the system will eventually recover.
The finite-time basin of attraction
\begin{equation}
    \B^T_\varepsilon(\A)=\{\x_0\in\B(\A)\,|\,\tau^{(\varepsilon)}(\x_0)\leq T\}
\end{equation}
includes all points in the basin of attraction for which the convergence time is short. The time horizon $T$ is the cutoff time after which a trajectory that is not yet $\varepsilon$-close to an attractor is not considered to be in any finite-time basin. The finite-time basin stability is simply defined as the weight of this set under uncertainty
\begin{equation}
    S^{(\nu,\varepsilon,T)}_{\mathrm{ft}} = \nu(\B^T_\varepsilon(\A)).
\end{equation}

% Add these in the preamble or just before the table
\renewcommand{\arraystretch}{1.5}  % Vertical padding
\setlength{\tabcolsep}{12pt}       % Horizontal padding

\begin{table*}[]
\begin{tabular}{lllp{6cm}}
\hline
\multicolumn{1}{|l|}{Resilience Measure} & \multicolumn{1}{l|}{Symbol} & \multicolumn{1}{l|}{References} & \multicolumn{1}{l|}{Description}                                                                                                                                                                 \\ \hline\hline
\multicolumn{4}{l}{Local Resilience}                                                                                                                                                                                                                                                                        \\ \hline
\hyperref[mes: crt]{Characteristic Return Time}               & $t_R$                       & \onlinecite{May1973StabilityEcology, Pimm1984Resilience, Holling1996EngineeringRV, Scheffer2009EWSNature, Arnoldi2016ResilienceMathComparison, Smith2022VegetationResilience}                                & Slowest asymptotic convergence rate in the linear regime of a point attractor. Computed by considering the largest real part of any eigenvalue of the Jacobian matrix: $-1/\lambda_\mathrm{max}$. \\
\hyperref[mes: reac]{Reactivity}                              & $R_0$                       & \onlinecite{Neubert1997AlternativeResilience, Verdy2008Reactivity, Snyder2010Reactivity, Arnoldi2016ResilienceMathComparison, Liu2021Reactivity, Capdevila2020Resilience}                                & Largest immediate response of the linear system around a point attractor to a unit disturbance in terms of the derivative of the absolute value at that disturbance.                    \\
\hyperref[mes: ampl]{Maximal Amplification (Time)}                    & $\rho_\mathrm{max},\,t_\mathrm{max}$         &  \onlinecite{Neubert1997AlternativeResilience, Arnoldi2016ResilienceMathComparison, Liu2021Reactivity, Capdevila2020Resilience}                               & Considering all trajectories started at unit distance from a point attractor in the linear regime, $\rho_\mathrm{max}$ gives the largest excursion away from the point attractor, while $t_\mathrm{max}$ gives the time at which it is attained.                                                                                                                         \\ \hline
\multicolumn{4}{l}{Nonlocal Resilience}                                                                                                                                                                                                                                                                    \\ \hline
\multicolumn{4}{l}{\textit{Geometrical measures}}                                                                                                                                                                                                                                                           \\
\hyperref[mes: mcs]{Minimal Critical Shock}                   & $s_\mathrm{min}$            & \onlinecite{Peterson1998Resilience, Beisner2003ResilienceLakes, Kerswell2014ResilienceOptApproach, Klinshov2015DistanceThresholdNumerical, Ashwin2016NetworkThresholds, Klinshov2020SwitchingThreshold, Halekotte2020MinimalFatalShocks}                                & Smallest distance of the attractor to any of the points on the boundary of its basin of attraction.                                                                                              \\
\hyperref[mes: mncs]{Maximal Non-Critical Shock}               & $s_\mathrm{max}$            &  introduced in this work                               & Largest distance of the attractor to a point of its own basin of attraction.                                                                                                                     \\
\hyperref[mes: bs]{Basin Stability}                          & $S$                         & \onlinecite{Holling1973ResilienceStabilityEcological, Gruemm1976DefinitionsResilience, Menck2013BasinStability, Lundstrom2018FindNonlocalResilience}                                & Volume of the basin of attraction of an attractor computed under a given probability measure of initial conditions.                                                                              \\

\multicolumn{4}{l}{\textit{Transient measures}}                                                                                                                                                                                                                                                             \\
\hyperref[mes: ct]{Convergence Time}                  & $\tau$                      &  \onlinecite{ONeill1976ResilienceEcosystem, DeAngelis1989ResilienceNutrients, Neubert1997AlternativeResilience, Cottingham1994Resilience}                               & Time it takes for an initial condition of positive probability to converge to an  $\varepsilon$-neighbourhood of the attractor. This information could be aggregated for many initial conditions by considering the mean, the maximum, or the median.                                                          \\
\hyperref[mes: cp]{Convergence Pace}                  & $\beta$                     &  introduced in this work                               & Equivalent to $\tau$ but considering instead the convergence time divided by the distance of the initial condition to the attractor. \\
\hyperref[mes: ftbs]{Finite-Time Basin Stability}              & $S_T$                       &   \onlinecite{Lundstrom2018FindNonlocalResilience, Schultz2018ShockFirstExitTime}                              & Equivalent to $S$ but considering only the finite-time basin, which contains initial conditions that converge to an $\varepsilon$-neighbourhood of the attractor in a finite-time $T$.         
                                     
\end{tabular}
\caption{List of all resilience measures for dynamical systems considered in this work. While local resilience measures are computed based on the linearisation of the system, nonlocal resilience measures take global state space and transient behaviour into account. Convergence time and pace are further synthesized into a basin median before any analysis.}\label{tab: measures}
\end{table*}

\section{Numerical Estimation} \label{sec: numEst}
The introduced resilience measures are mathematically well-defined and assume values in the real numbers or infinity. However, they are generally not computable analytically. Even finding point attractors $\A=\{\x^*\}$, i.e., solving for $f(\x^*,p)=\mathbf{0}$ and checking for linear stability, is only possible analytically in the simplest of cases. 
For this reason, a numerical treatment of the dynamical system is necessary.

In this work, we develop a new numerical algorithm that can estimate all resilience measures listed in the previous section.
The estimation occurs in parallel, accumulating only the minimum amount of information relevant to this goal. As a result, the algorithm is not only efficient but also easily extendable due to its design. New measures can be added to the routine as they arise.
We have implemented this algorithm as an open source contribution to the DynamicalSystems.jl library, \cite{Datseris2018DynamicalSystemsjl} in its submodule Attractors.jl. \cite{Datseris2023FrameworkGlobalStability} This submodule also contains the routines we employ to find attractors and their basins, using the accumulating recurrences method of Ref.~\onlinecite{Datseris2022BasinsAttraction}. Specifically, this infrastructure allows mapping individual initial conditions to new or already identified attractors by terminating trajectories only after a specified number of recurrences on a grid have taken place. This grid and the initial conditions employed for finding all attractors have to be chosen with contextual knowledge of the system. They are, however, not essential to the estimation of resilience measures.

The estimation of resilience measures takes part in a second step, using the identified attractors as a starting point. For the set of local resilience measures, only point attractors are considered. The associated Jacobian matrix is the basis of all local measures and is estimated using the Julia package ForwardDiff.jl. From there, the characteristic return time can be estimated as the spectral abscissa of the eigenvalues of the Jacobian, to give one example. The computation of reactivity and maximal amplification (time) similarly follows their defining equations in Sec.~\ref{sec: resMes}.

\begin{figure}[t]
	\centering
	\includegraphics[width=\columnwidth]{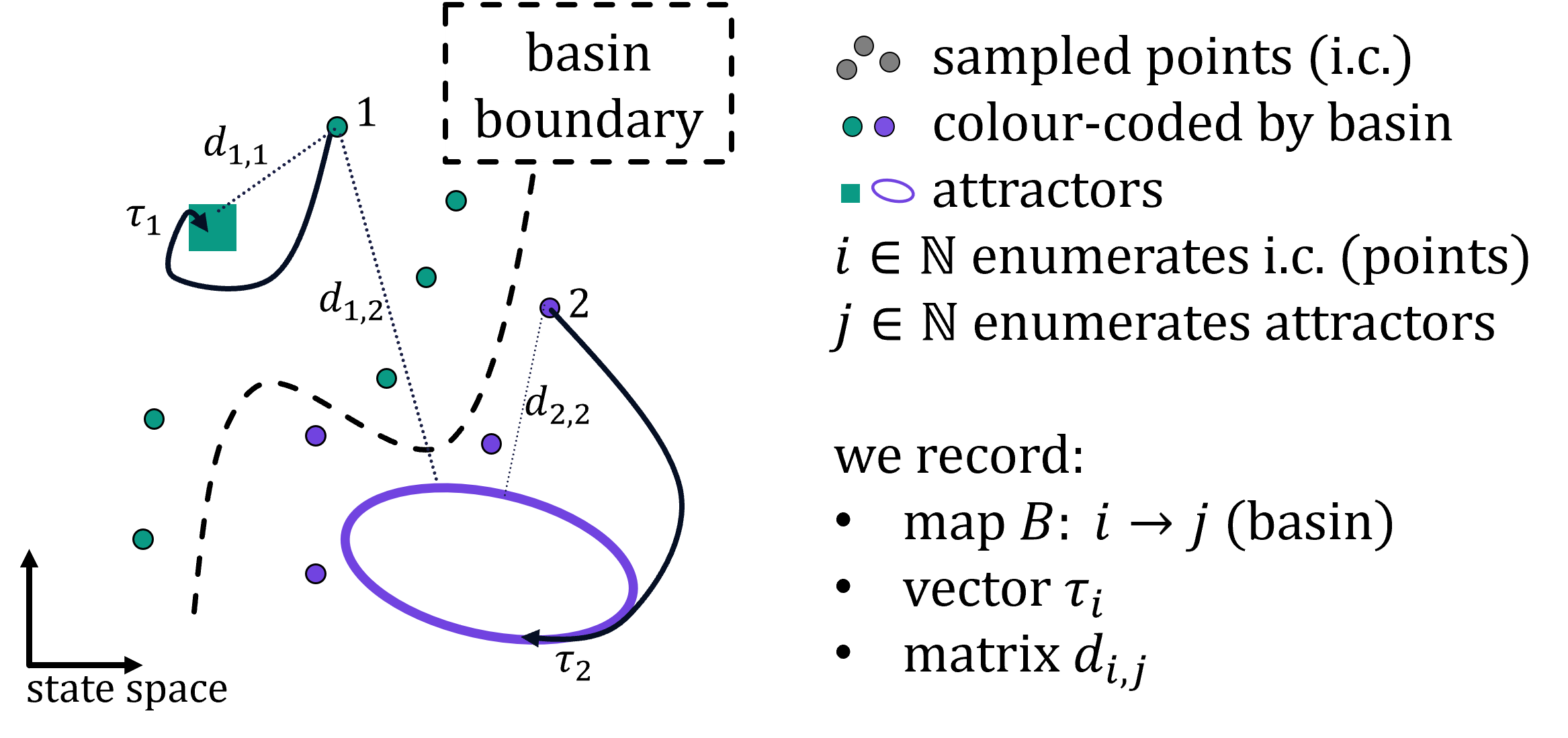}
	\caption{Schematic illustrating the key quantities recorded to calculate different measures of resilience for a given dynamical system. The information we record while mapping initial conditions (enumerated by $i\leq N$) to attractors (enumerated by $j \in \mathbb{N}$, with $j = 0$ meaning divergence) are: (1) The initial conditions $\mathbf{x}_{i,0}$, and the basin mapping $B: i \to j$, stored efficiently as a vector, using the fact that $i$ runs from 1 to $N$.
    (2) The convergence time of each initial condition $\tau_i$, stored as a vector. In our applications, this is defined as the first time where the trajectory is $\varepsilon$-close to the attractor $\A_j$. (3) A matrix of distances between initial conditions $\x_{i,0}$ and attractors $\A_j$. Since we are using the Euclidean norm in the ordinary distance between sets, the matrix stores $d_{i,j}=\min_{\y\in\A_j}\norm{\x_{i,0}-\y}_2$. However, the algorithm allows for arbitrary choices of distances.
    }
    \label{fig: algorithm}
\end{figure}

For the nonlocal measures of resilience, the goal of the remaining procedure is to gather information needed for their computation efficiently and in parallel. 
This is done by again mapping initial conditions to attractors while recording some additional information. 
In the software implementation, this mapping may continue using the same recurrence-based algorithm.
Alternatively, which is also the approach we use in this work, this mapping routine can be changed to a different termination condition. Trajectories are considered to have converged to an attractor not when a set number of recurrences has taken place, but when $\varepsilon$-proximity to an attractor has been reached. 
We thus gather information about which basin of attraction an initial condition belongs to, while obtaining the convergence time in the sense of Eq.~\eqref{eq: ct} at the same time. 
A third information continuously stored along the routine is the distance of an initial condition to its attractor, which is needed for the computation of the convergence pace and critical shocks.
We portray the three pieces of information in Fig.~\ref{fig: algorithm}, which illustrates the computational setup and details the computed values.

From here, the computation of the resilience measures on a finite amount of data from initial conditions is straightforward.
To give just two examples, one can estimate the basin stability, i.e., the distribution weighted volume of basins of attraction, from $N$ initial conditions by calculating
\begin{equation}
    S(\A_j) \approx \frac{1}{N}\sum_i w_\nu(\mathbf{x}_{i,0})\delta_{B(i), j}
\end{equation}
with $\delta_{i,j}$ the Kronecker delta. The minimal critical shock can be computed as
\begin{equation}
s_{\mathrm{min}}(\A_j) \approx \min_{1\leq i\leq N}(d(i,j)| B(i) \ne j).
\end{equation}

In fact, any resilience measure that can be estimated through the three recorded quantities $B, \tau, d$ can be directly integrated into our algorithm and be estimated essentially for free.

The number of initial conditions $N$ is a central configuration specifier. In the routine, $N$ points are uniformly sampled at random from a set state space region of interest. The weighting distribution $\nu$ only enters the computation of resilience measures at a later point. For a thorough exploration of the basins of attraction via such sampling, a sufficient amount is necessary. With increasing dimensionality, this condition is likely to increase the computational cost of the routine drastically. In Appendix \ref{app: synthetic}, we investigate the dependence of the accuracy of estimated resilience measures on $N$ for a problem of fixed dimensionality $2$ to give an initial impression of the needed number of initial conditions.

The process outlined above is for a particular parameter value $p$, determining the ODE right-hand side $f(\cdot,p)$. When combined with the global continuation, the resilience measures per attractor are matched to their previous values using the same matching technique outlined in Ref.~\onlinecite{Datseris2023FrameworkGlobalStability}. This allows for a consistent tracking of resilience as a system moves towards a qualitative change via parameter variation, such as a bifurcation.

The code implementation is straightforward to use. In Appendix~\ref{app: code} we present a code snippet that performs all numerics discussed in this section.

\section{Results} \label{sec: results}

We employ the numerical estimators of the introduced resilience measures to three systems from the applied sciences. These will be a box model of the Atlantic Meridional Overturning Circulation (AMOC), \cite{Alkhayuon2019AMOCRateTipping, Wood20195BoxModel} a predator-prey population model, \cite{Zeng2024PredPreyHollingIII} and the Lorenz-84 model of atmospheric turbulence. \cite{Lorenz1984Model} In each of these cases, we determine one system parameter of interest in the context of potentially catastrophic regime shifts. 
We note that the global continuation provided by Attractors.jl would also allow for any arbitrary parameter curve in parameter space to be tracked instead. 
The changes to the system's dynamics in terms of its evolving resilience measures are then assessed, providing new perspectives on these well-established conceptual models.

For all three examples, we choose a probability distribution of initial conditions that is uniform on a rectangular state-space region (cuboid for the three-dimensional example). These are the same regions illustrated in Figs.~\ref{fig: AMOCBasin} and \ref{fig: predpreyBasin}, where the attractors and basins of attraction are visualized. For the third example, the state-space is given in the description of Attractor visualization Fig.~\ref{fig: lorenz84Attractors}. In each example, we chose to sample $N=10^5$ initial conditions in the respective state space regions for every parameter value $p$. The convergence threshold $\varepsilon$ and finite-time horizon $T$ are system-specific and defined during the analyses below.

\subsection{AMOC box model}
The Atlantic Meridional Overturning Circulation (AMOC) is largely responsible for the transport of warm, saline water from the Tropics to the North Atlantic, thus playing an essential role in Europe's climate. \cite{Johns2023AMOCHeatTransport} The AMOC has been suspected to be influenced by an increased freshwater flux in the North Atlantic. \cite{Rahmstorf2005AMOCHystersis, VanWesten2023AMOCHysteresis} This is currently brought about by climate change and the associated Greenland Ice Sheet melt. \cite{Yang2016AMOCFreshwater} Simplified box models represent the dynamics of the circulation in a coarse partitioning of the world's oceans. \cite{Stommel1961CircModel, Wood20195BoxModel} Its parameters can be tuned so that the model reproduces the dynamics observed in more comprehensive Earth System models as best as possible. \cite{Chapman2024AMOCNoiseParameters, Alkhayuon2019AMOCRateTipping} Investigating the box model's dynamical characteristics with respect to a larger freshwater influx is thus conceptually valuable in the analysis of the real-world system. We track the evolution of resilience measures in this context and thereby provide insights about the system's response to climate change. 

The 3-box model investigated here was introduced in Ref.~\onlinecite{Alkhayuon2019AMOCRateTipping} as a simplification of the 5-box model first proposed in Ref.~\onlinecite{Wood20195BoxModel}. It consists of a Northern, Tropical, and Southern Atlantic box, and the concentrations of salt are described by a two-dimensional ODE. This is because salt is conserved over all three boxes, so that the system may be reduced by one dimension.
\begin{widetext}
    \begin{align}
    V_N \frac{\dint S_N}{\dint t}&=q\left[\left(\frac{K_N}{q}+\Theta(q)\right)(S_T-S_N)-\Theta(-q)(S_B-S_N)\right]-(F_N+aH)S_0\nonumber\\
    V_T \frac{\dint S_T}{\dint t}&=q\left[\Theta(q)(\gamma S_S+(1-\gamma)S_{IP}-S_T)-\left(\Theta(-q)-\frac{K_N}{q}\right)(S_N-S_T)\right]+K_S(S_S-S_T)-(F_T+bH)S_0\label{eq: AMOC}\\
    q&=\frac{\lambda(\alpha(T_S-T_0)+\beta(S_N-S_S))}{1+\lambda\alpha\mu}\nonumber
    \end{align}
\end{widetext}
Here, $S_N$ and $S_T$ are the dynamical salinity quantities in the Northern and Tropical Atlantic boxes, $q$ describes the flow strength of the AMOC, and $\Theta$ is the Heaviside function. The parameter $H$ represents the freshwater hosing induced by climate change, which is suspected to have potentially catastrophic consequences at a certain threshold, or tipping point. All other quantities are fixed parameters.

\begin{figure}[b]
	\centering
	\includegraphics[width=\columnwidth]{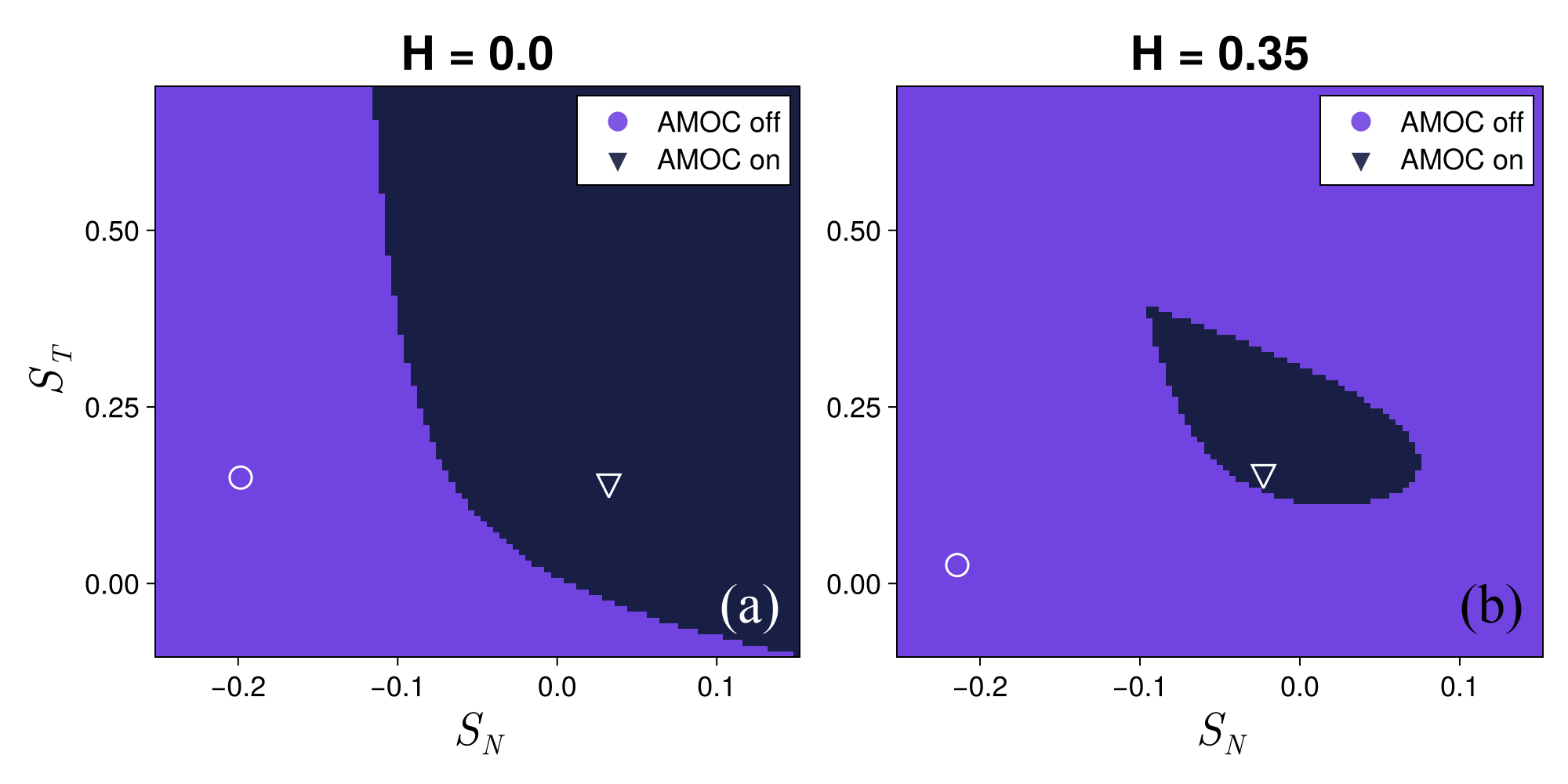}
	\caption{Attractors and basins of attraction for the AMOC 3-box model of Eq.~\eqref{eq: AMOC}. The two-dimensional system is analysed for two values of the freshwater hosing parameter $H$. Values of all other fixed parameters can be found in the code documentation. Under stronger hosing conditions, the AMOC on state loses resilience in the sense of its basin stability. We quantify this loss and the evolution of other resilience measures numerically in Fig.~\ref{fig: AMOC}.}
    \label{fig: AMOCBasin}
\end{figure}

Studying the characteristics of this system in its parametrization corresponding to zero freshwater hosing, two stable states of the AMOC can be found. One corresponds to the currently observed strong AMOC (on), while the other represents a weak or even reversed AMOC (off). We can visualise this in the 2D plane of dynamic salinities by numerically generating the basins of attraction belonging to these attractors (see Fig.~\ref{fig: AMOCBasin}(a)). Initial conditions in this plane represent certain salinities in the Northern and Tropical Atlantic and will generally either converge to the AMOC on state or the off state. In case they lie on the zero-volume boundary of the two cases, they may linger there for an infinitely long time. \cite{Lucarini2017EdgeTracking, Lohmann2025EdgeStatesEWS} Generating the same figure for a parameter setting with a stronger North Atlantic fresh water flux, one can observe a change in the basin structure (see Fig.~\ref{fig: AMOCBasin}(b)). By the naked eye, arguably, the on state attractor has lost resilience, since fewer initial conditions converge to it.

\begin{figure}[H]
	\centering
	\includegraphics[width=\columnwidth]{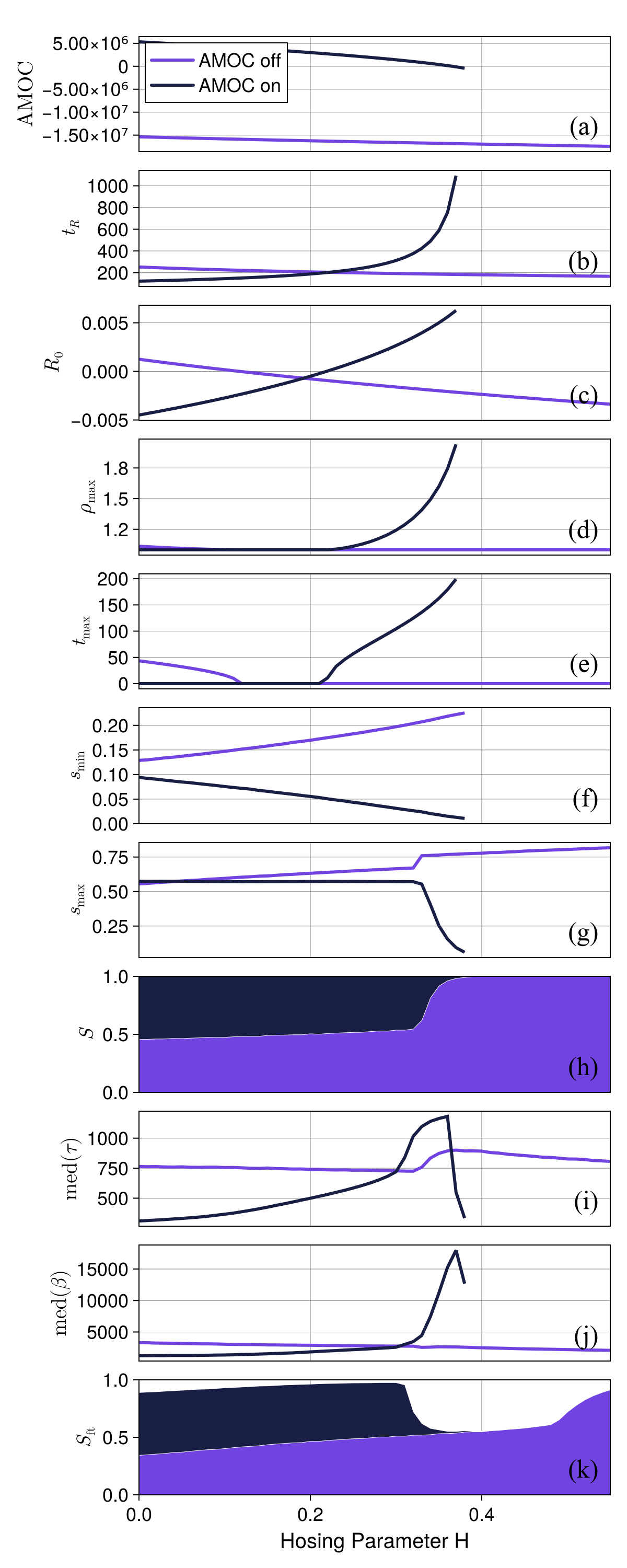}
	\caption{Resilience measures related to the on and off state of the AMOC 3-box model from Eq.~\eqref{eq: AMOC}. In panel (a), the AMOC flow strength $q$ is illustrated. Color scheme identical to Fig.~\ref{fig: AMOCBasin}.}
    \label{fig: AMOC}
\end{figure}

However, this is just one of the interpretations of system resilience (basin stability), and we have set out in this work to quantitatively investigate the evolution of a broader class of resilience measures. Figure \ref{fig: AMOC} shows the numerical estimation of the resilience measures introduced in Section \ref{sec: resMes}. Their respective trends are most interesting in advance of the vanishing of the stable on state by means of a Hopf bifurcation at the freshwater parameter value $H_\mathrm{Hopf}\approx0.39$. The system undergoes a non-catastrophic global bifurcation shortly before this catastrophic one: At $H_\mathrm{hom}\approx0.36$, an unstable limit cycle emerges from a homoclinic bifurcation at the basin boundary of the on state. This limit cycle vanishes again at the subcritical Hopf bifurcation. As expected for a Hopf bifurcation, the size of the basin of attraction, as measured by the basin stability, shrinks to 0 when the AMOC on state vanishes (see Fig.~\ref{fig: AMOC}(h)). Initial conditions also seem to converge more slowly, as can be seen by the evolution of the median convergence time and pace, though we observe a precipitous fall immediately before the bifurcation (see Fig.~\ref{fig: AMOC}(i, j)). This could be due to the fact that as the basin closes in on the on state, initial conditions with long transients are lost, decreasing the median time and pace. Using a convergence threshold of $\varepsilon=0.01$ and a finite-time horizon of $T=1000$, the comparison of basin stability and finite-time basin stability reveals a qualitative difference in the longest convergence times of the two basins. More specifically, after the collapse of the AMOC on state, the state space newly allocated to the off state's basin exhibits a prohibitively long convergence time, which slowly adjusts for larger parameter values, raising the finite-time basin stability (see Fig.~\ref{fig: AMOC}(k)). In terms of local resilience, the on state unsurprisingly becomes linearly unstable at the bifurcation point ($t_R\rightarrow\infty$ in Fig.~\ref{fig: AMOC}(b)). The finite-time responses of the linear system also change, with the reactivity rising smoothly while the amplification logically commences its rise only at the time where the reactivity crosses $0$ (see Fig.~\ref{fig: AMOC}(c--e)). This crossing point occurs at an inconspicuous time, considering the global qualitative changes and the associated bifurcation thresholds.

\newpage%%%%%%%%%%%%%%%%%%%%%%%%%%%%%%%%%%
\subsection{Predator-prey model}
We now consider another two-dimensional dynamical system and estimate the resilience of its attractors. The general structure of a predator-prey population model of Holling type III with Allee effect is given by
\begin{align}
    \frac{\dint x}{\dint t}&=x(1-x)(x-E)-\frac{x^2y}{Ax^2+Bx+1},\\
    \frac{\dint y}{\dint t}&=C\frac{x^2y}{Ax^2+Bx+1}-Dy\label{eq: predprey}.
\end{align}
The dynamics of the prey population size $x$ and the predator population size $y$ are constrained by a growth and death rate, respectively. These are given by the functional relationships $x(1-x)(x-E)$ and $-Dy$. Additionally, the consumption of the prey by the predator takes place at a rate $x^2y/(Ax^2+Bx+1)$, which enhances the predator population by the same rate times a factor $C$. The behaviour of this dynamical system has been extensively investigated with respect to each of the five parameters. \cite{Zeng2024PredPreyHollingIII} We are interested in the effects of a change of the Allee growth rate $x(1-x)(x-E)$ via the parameter $E$.

\begin{figure}[b]
	\centering
	\includegraphics[width=\columnwidth]{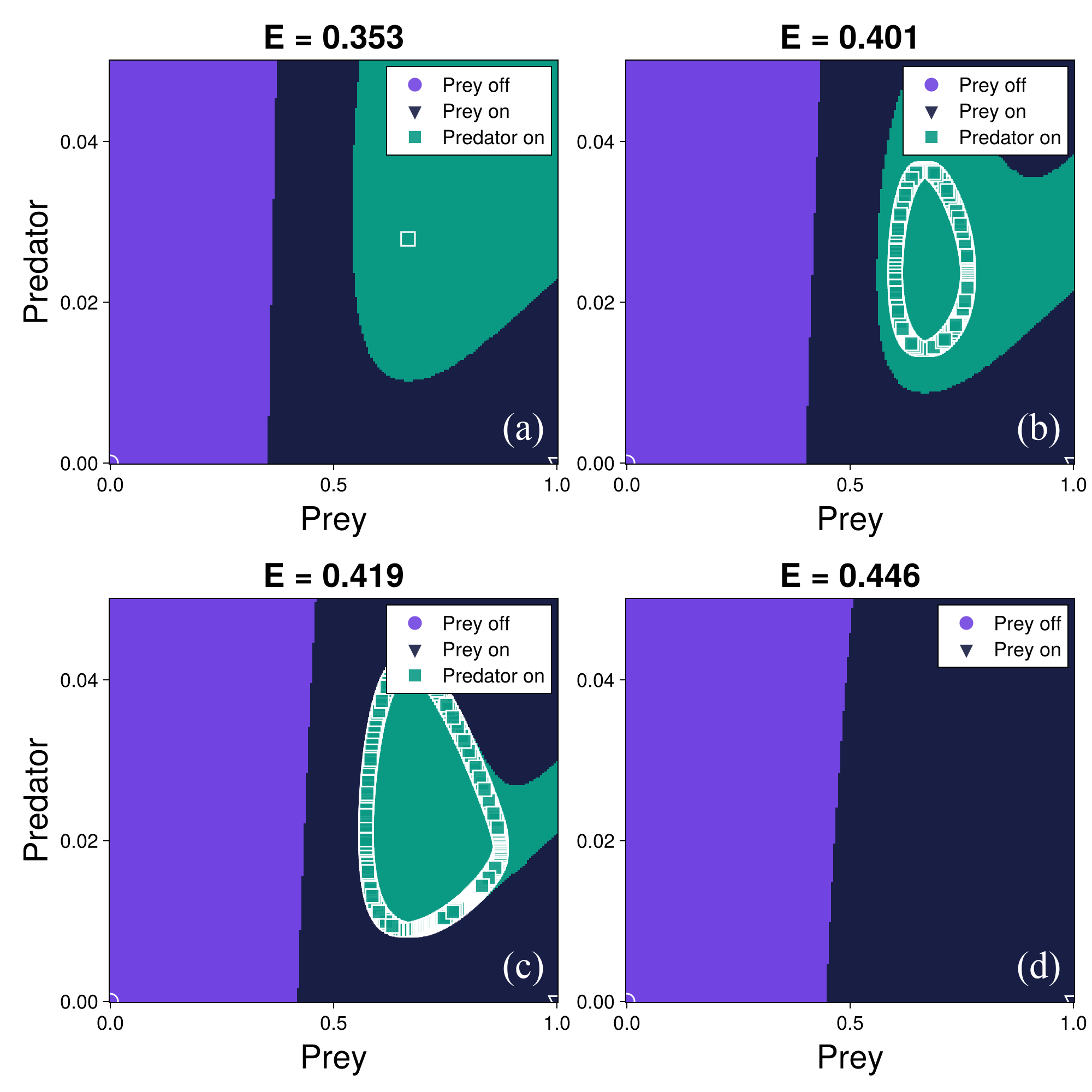}
	\caption{Attractors and basins of attraction for the predator-prey model introduced in Eq.~\eqref{eq: predprey}. Parameter values are set to $A=2.05,\,B=-2.6,\,C=0.4,\,D=1$. The two-dimensional system is analysed for four values of the Allee effect parameter $E$. The desired state of a functioning ecosystem with active predator and prey populations undergoes a Hopf bifurcation, turning the point attractor into a limit cycle. For even larger values of $E$, this limit cycle and the associated system functioning vanish entirely. This occurs through a merging of the limit cycle with the basin boundary. We track the numerically estimated resilience measures through these structural changes in Fig.~\ref{fig: predprey}.}
    \label{fig: predpreyBasin}
\end{figure}

\begin{figure}[H]
	\centering
	\includegraphics[width=\columnwidth]{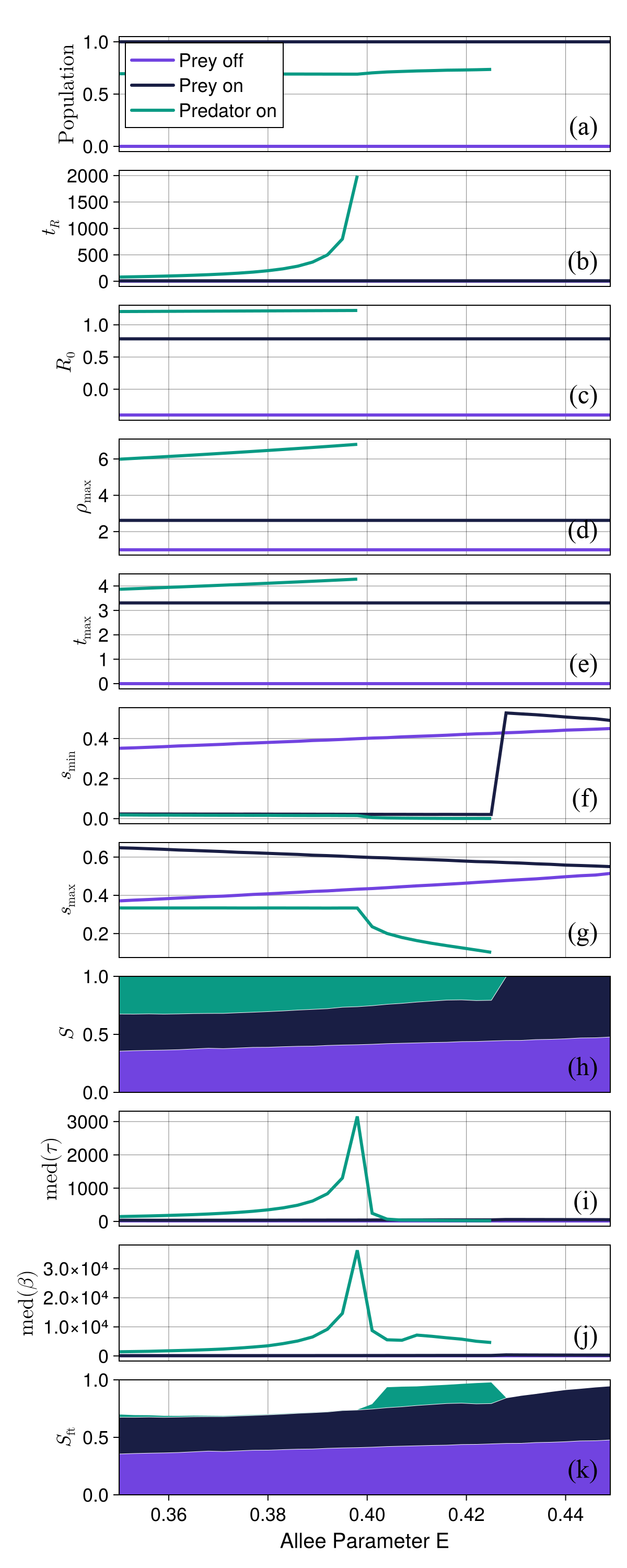}
	\caption{Resilience measures related to the three stable states of the predator-prey model introduced in Eq.~\eqref{eq: predprey}. Panel (a) shows the total sum of population sizes. Color scheme identical to Fig.~\ref{fig: predpreyBasin}.}
    \label{fig: predprey}
\end{figure}

Figure \ref{fig: predpreyBasin} shows the basin configurations for four exemplary values of $E$. 
In all cases, there exist two stable point attractors at the edges of the considered population domain. These represent total extinction at $(x,y)=(0,0)$, and success of the prey but extinction of the predator at $(1,0)$. Additionally, there exists a regime of non-zero population sizes. It is this non-trivial state that we would like to consider with respect to its resilience. For low values of $E$, this state is a point attractor, which transitions to a stable limit cycle in a supercritical Hopf bifurcation at $E_\mathrm{Hopf}\approx0.398$. At a critical value of $E_c\approx0.423$, this mode of system functioning ceases to exist because of a boundary crisis in the basin of attraction. That is, the stable orbit merges with the basin boundary and vanishes. For high values of $E$, the system thus only exhibits the discussed trivial states without species cohabitation. Through the course of these occurrences, we track the estimated resilience measures of all attractors. For non-point attractors, the local resilience measures introduced in this work are not applicable. Figure \ref{fig: predprey} shows the estimated nonlocal resilience measures and the local measures where applicable.

The system state is quantified as the mean population sum of both species on the respective attractors. For the two trivial point attractors, this population sum does not change while the non-trivial state, transformed through the Hopf bifurcation, changes slightly. At the bifurcation, we again see a loss of linear stability for the state of interest. Reactivity appears to be constant above $1$, while maximal amplification (time) rises steadily. The investigation of convergence times through their median values and the finite-time basin stability again gives an insight into the changes to the transient behaviour of the system. Here, we chose $\varepsilon=0.001$ and $T=100$. 
The median convergence time and pace rise monotonically before the system undergoes a Hopf bifurcation and abruptly decrease after. Interestingly, the decrease in the convergence pace appears less precipitous than that of the convergence time. This drop seems to come with a substantial part of the state space now exhibiting convergence times below the finite-time horizon. Finite-time basin stability is non-zero in this parameter range between the Hopf bifurcation and the eventual boundary crisis. One more qualitative difference observed in the vanishing of the periodic orbit compared to that of the AMOC on state of the last example can be seen when investigating the basin stability. Instead of a continuous decrease towards 0, the basin stability suddenly drops from a finite value when the boundary crisis occurs. This particular resilience measure is thus not a reliable indicator of an impending collapse in this application.

\newpage%%%%%%%%%%%%%%%%%%%%%%%%%%%%%%%%%%%%
\subsection{Lorenz-84 model}
\begin{figure}[b]
	\centering
	\includegraphics[width=\columnwidth]{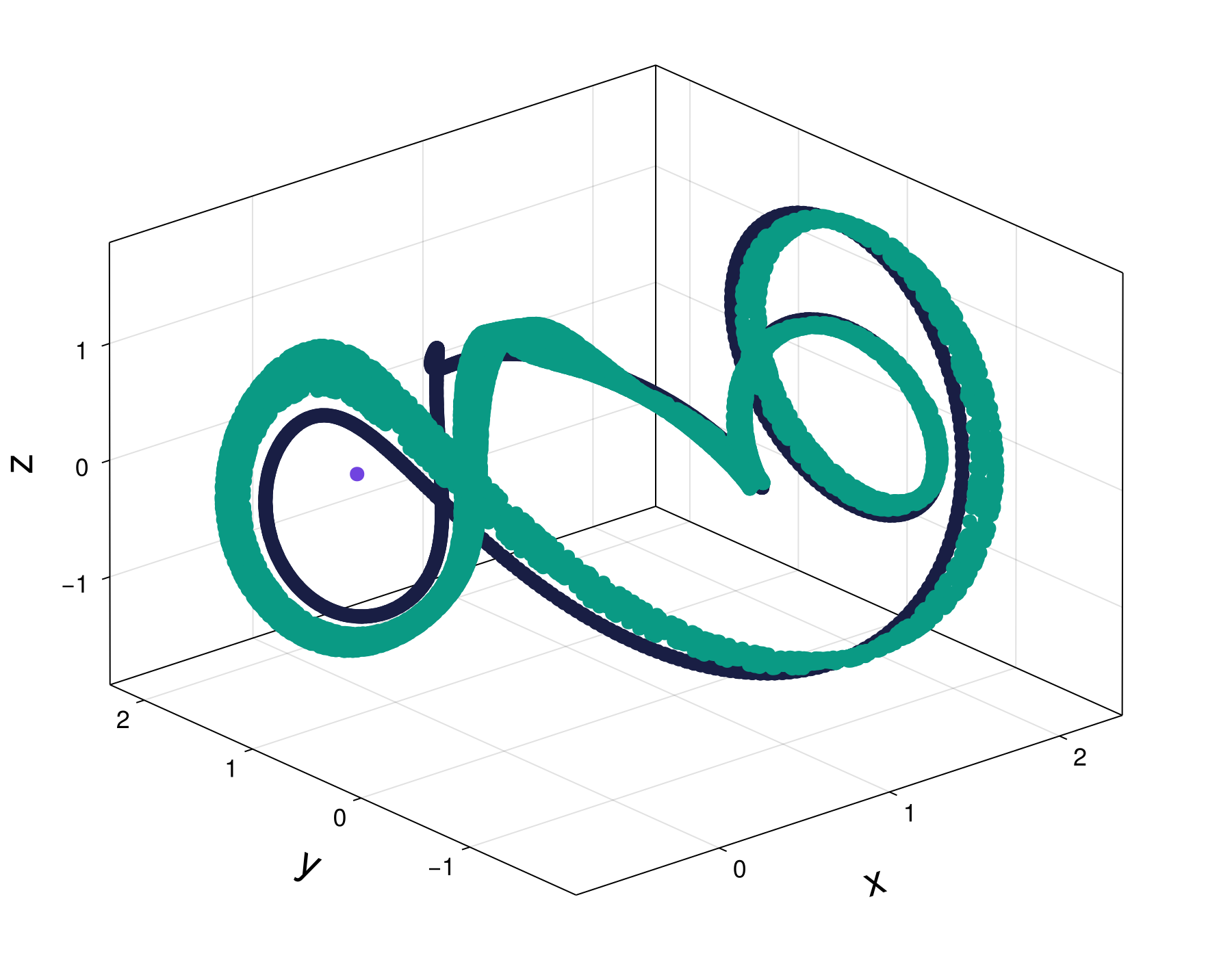}
	\caption{Attractors of the Lorenz-84 system defined in Eq.~\eqref{eq: Lorenz84} with parameter values $F=6.886,\,G=1.355,\,a=0.255,\,b=4.0$. At this particular value of $G$, the system has three attractors. One point attractor, one limit cycle, and one chaotic attractor. The latter two will merge into a single limit cycle for larger values of $G$. We vary this parameter and estimate the resilience of all attractors of the system in Fig.~\ref{fig: lorenz84}. The considered state space region for sampling the initial conditions is delimited by $-4$ and $4$ in the $x$, $y$, and $z$ dimension, respectively.}
    \label{fig: lorenz84Attractors}
\end{figure}

As a third example, we investigate a system that exhibits a transition to chaos. The Lorenz-84 model is a simplified low-order model of atmospheric circulation. \cite{Lorenz1984Model} It was designed to capture essential features of mid-latitude weather dynamics, particularly the interaction between large-scale westerly winds and smaller-scale traveling waves in the atmosphere. Its evolution equations read
\begin{align}
\frac{\dint x}{\dint t} &= -y^2 - z^2 - ax + aF, \\
\frac{\dint y}{\dint t} &= xy - bxz - y + G, \label{eq: Lorenz84}\\
\frac{\dint z}{\dint t} &= bxy + xz - z.
\end{align}
The system consists of three coupled nonlinear differential equations representing the mean westerly wind $x$, and two components of planetary-scale wave activity $y$ and $z$. Despite its simplicity, the model displays a rich variety of dynamical behaviors, including equilibrium points, periodic oscillations, and chaos \cite{Shilnikov1995Lorenz84Bifurcation} (see Fig.~\ref{fig: lorenz84Attractors}).

A key value within the Lorenz-84 model is the parameter $G$, which represents the invariant wave forcing due to eddy activity at unresolved scales. As $G$ changes, the system undergoes several bifurcations that lead to dynamics of differing complexity and, at times, chaotic behavior. 
Transitions to chaos can be interpreted as the emergence of turbulence in the atmosphere, which marks an important qualitative change. 
\begin{figure}[H]
	\centering
	\includegraphics[width=\columnwidth]{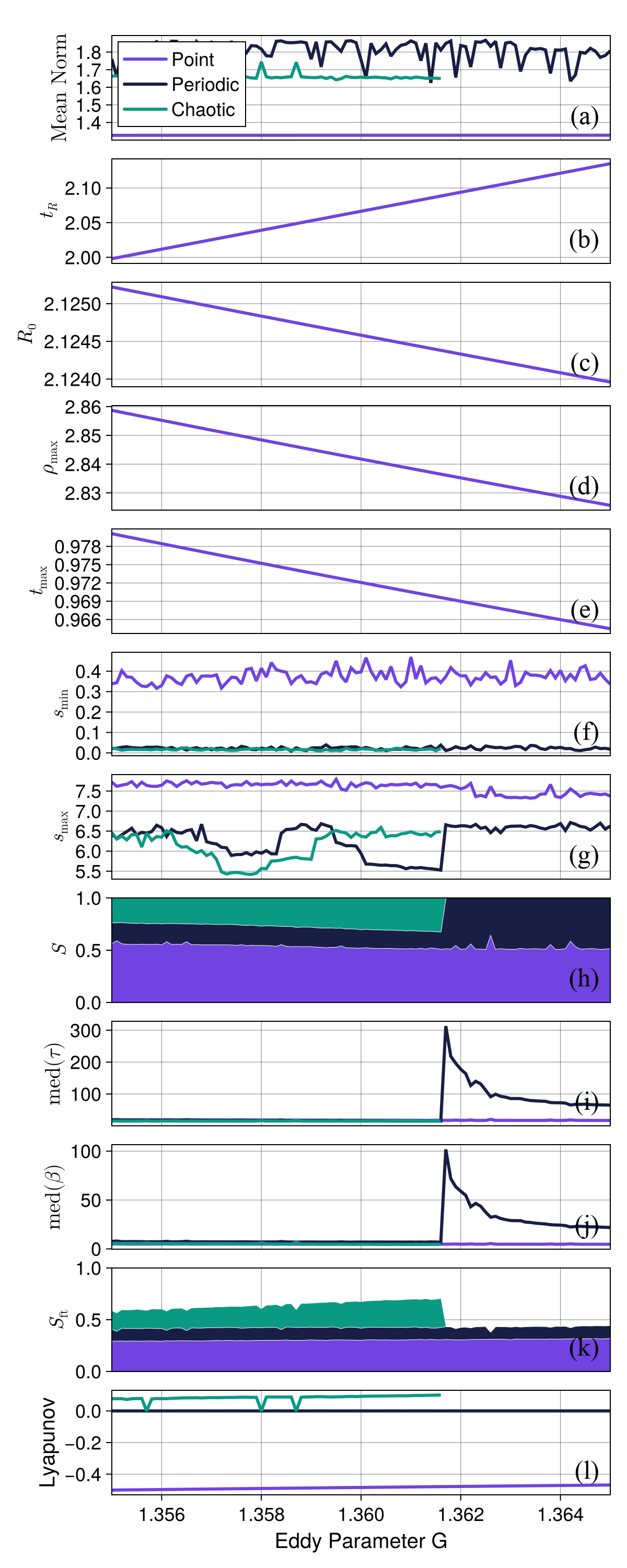}
	\caption{Resilience measures related to the Lorenz-84 system introduced in Eq.~\eqref{eq: Lorenz84}. The system exhibits up to three attractors in the considered parameter range, which are represented in panel (a) by the average value of the norm of constituent points. Color scheme identical to Fig.~\ref{fig: lorenz84Attractors}.}
    \label{fig: lorenz84}
\end{figure}
We track the evolution of our selection of resilience measures as the system goes through one such change. Figure \ref{fig: lorenz84} shows the results. The three different attractors of the system seen at the initial parameter value in Fig.~\ref{fig: lorenz84Attractors} are of different natures. One of the attractors is a stable equilibrium point, while the others are periodic and chaotic attractors, respectively. The chaotic attractor vanishes at the parameter value $G_c\approx1.6318$ by merging with the periodic orbit. It is therefore also of interest to investigate the behaviour of the system on the right-hand side of this parameter value, as this precedes the emergence of chaos when reducing $G$. We have chosen to additionally show the estimated maximal Lyapunov exponent of each attractor to differentiate the attractors' qualities (see Fig.~\ref{fig: lorenz84}(l)). It can be seen that at several parameter values, the maximal Lyapunov exponent of the chaotic attractor drops to $0$. Upon closer inspection, this does not appear to be an estimation inaccuracy but to correspond to short crises and microscopic periodic windows of that attractor.

The changes to the linearised dynamics around the point attractor are monotonous as a function of the parameter and appear unaffected by the qualitative change of the attractor structure at $G_c$. Due to the more convoluted attractor structures in three-dimensional space, the basin geometry indicators of minimal and maximal critical shock exhibit more complex evolutions than in the preceding examples. 
The median convergence time and pace of the periodic orbit increase in advance of the emergence of the chaotic attractor, i.e., when decreasing $G$ towards $G_c$. This can be interpreted as a decrease in system resilience, as trajectories exhibit longer transients before returning to a familiar mode of functioning. With convergence threshold $\varepsilon=0.01$ and finite-time horizon $T=20$, we see a non-proportional scaling when moving from the basin stability to the finite-time basin stability. While nearly all of the basin of the chaotic attractor qualifies as its finite-time basin, the basins of the other two attractors are substantially diminished when imposing a cut-off convergence time. Interestingly, all of the basin volume that belonged to the chaotic attractor at the time of its disappearance seems henceforth to be associated with long transients that eventually converge to the joint periodic attractor.

\section{Conclusions and Outlook}
\renewcommand{\arraystretch}{1.5}  % Vertical padding
\setlength{\tabcolsep}{12pt}
\begin{table*}[]
\begin{tabular}{lp{6cm}ccccc}
\hline
\multicolumn{1}{|l|}{Measure}                             & \multicolumn{1}{l|}{Limitations}                                                                                                                                                                       & \multicolumn{5}{l|}{Observed behaviour before regime change}                             \\ \hline
Local                                                     &                                                                                                                                                                                                        & AMOC          & \multicolumn{2}{c}{Predator-Prey} & \multicolumn{2}{c}{Lorenz-84}        \\ \hline
                                                          &                                                                                                                                                                                                        & \textit{Hopf} & \textit{Hopf}   & \textit{Hom.}   & \textit{Distruct.} & \textit{Emerg.} \\
\hyperref[mes: crt]{$t_R$}                                & Only defined for point attractors. First-order smoothness of ODE required. Estimation reliant on the accuracy of ForwarDiff.jacobian functionality.                                                        & $\nearrow$    & $\nearrow$      & n.a.            & n.a.               & n.a.            \\
\hyperref[mes: reac]{$R_0$}                               & See above.                                                                                                                                                                                             & $\nearrow$    & $\rightarrow$   & n.a.            & n.a.               & n.a.            \\
\hyperref[mes: ampl]{$\rho_\mathrm{max},t_\mathrm{max}$} & Additionally to the above, complications can arise from the needed maximization routine. The global maximum in Eq.~\eqref{eq: amplification} will only be found with favorable initial tries of $t$. & $\nearrow$    & $\nearrow$      & n.a.            & n.a.               & n.a.            \\ \hline
Nonlocal                                                  &                                                                                                                                                                                                        & AMOC          & \multicolumn{2}{c}{Predator-Prey} & \multicolumn{2}{c}{Lorenz-84}        \\ \hline
\textit{Geometrical}                                      &                                                                                                                                                                                                        & \textit{Hopf} & \textit{Hopf}   & \textit{Hom.}   & \textit{Distruct.} & \textit{Emerg.} \\
\hyperref[mes: mcs]{$s_\mathrm{min}$}                     & Accurate exploration of the basin boundaries requires many sampled initial conditions, leading to a long runtime.                                                                                      & $\searrow$    & $\rightarrow$   & $\rightarrow$   & $\rightarrow$      & $\rightarrow$   \\
\hyperref[mes: mncs]{$s_\mathrm{max}$}                    & See above.                                                                                                                                                                                             & $\searrow$    & $\rightarrow$   & $\searrow$      & $\searrow$      & $\rightarrow$   \\
\hyperref[mes: bs]{$S$}                                   & Reliant on many sampled initial conditions but less than the above due to arguments outlined in Appendix \ref{app: synthetic}.                                                                         & $\searrow$    & $\searrow$      & $\searrow$      & $\nearrow$         & $\rightarrow$   \\
\textit{Transient}                                        &                                                                                                                                                                                                        & \textit{Hopf} & \textit{Hopf}   & \textit{Hom.}   & \textit{Distruct.} & \textit{Emerg.} \\
\hyperref[mes: ct]{$\mathrm{med}[\tau]$}                  & Even less reliant on the number of sampled initial conditions. For computing the mean or maximal value, a more thorough exploration of the tail of the distribution would be needed.                   & $\nearrow$    & $\nearrow$      & $\rightarrow$   & $\rightarrow$      & $\nearrow$      \\
\hyperref[mes: cp]{$\mathrm{med}[\beta]$}                 & See above.                                                                                                                                                                                             & $\nearrow$    & $\nearrow$      & $\searrow$      & $\rightarrow$      & $\nearrow$      \\
\hyperref[mes: ftbs]{$S_T$}                               & See $S$.                                                                                                                                                                                               & $\searrow$    & $\rightarrow$   & $\rightarrow$   & $\nearrow$         & $\rightarrow$  
\end{tabular}
\caption{Overview of the results derived from the application of our numerical routine to three example systems in Section \ref{sec: results}. In the second column, we report briefly on any limitations encountered during the estimation. In the third column, we describe the qualitative change of the measures in the advance of the respective regime shifts observed in the examples. These are: the Hopf bifurcation of the AMOC on state effectuating its disappearance, the Hopf and homoclinic bifurcations of the non-trivial population state in the predator prey model, and the destruction and emergence of chaos in the Lorenz84 model. For the last two transitions, we thus inspect the behaviour of the chaotic state before its disappearance (from left to right on the parameter axis) and that of the periodic state before the emergence of chaos (from right to left on the parameter axis). Arrows indicate the signs of the qualitative change of the respective measures.}\label{tab: effects}
\end{table*}
In this work, we have introduced and showcased a novel numerical framework for the estimation of resilience measures in dynamical systems. We have considered a relevant selection of such resilience measures from the literature and have discussed their heuristic interpretation. The open-source contribution to an existing and well-established Julia software package is convenient, efficient, and easily extended with new measures of resilience. For three example systems from applied contexts, we have tracked the estimated resilience of a desired mode of system functioning until its collapse in the form of a dynamic bifurcation. The behaviour of the system in advance of the catastrophic collapse or qualitative change can thereby be viewed through a novel and more comprehensive lens. We have summarized the qualitative results of our analyses in Table \ref{tab: effects}.

While we could not immediately extract any universal law for the asymptotics of resilience measures at a bifurcation point based on our computations, these considerations may spark future mathematical and applied research.

The numerical framework, though efficient in a relative sense, still comes with several practical drawbacks. Mapping out basins of attraction in state space inherently requires the exploration of all system dimensions as initial conditions. While this is feasible for systems of a handful of dimensions using the resource-efficient routines provided by Attractors.jl, the exponential increase in the number of initial conditions prohibits this for high-dimensional systems. For such applications, a prior dimension reduction could alleviate this problem. Short of such a drastic simplification, it may also be possible to identify less relevant dimensions of an anisotropic dynamical system. For instance, if some dimensions of the system are not relevant to the exit from any basin of attraction, they may be disregarded in the exploration of the state space for this purpose. \cite{Kuehn2015CurseInstability, Kuznetsov2021AttractorDimension, Kutz2016DynamicModeDecomposition} An intelligent Monte-Carlo sampling of initial conditions with respect to the given weighting distribution may ease the computation of resilience measures that are not based on the geometry of the basin of attraction but on mean values, such as basin stability. \cite{Hastings1970MonteCarlo, Ogata1989MonteCarloIntegration}

Another practical limitation arises from the a priori unconstrained estimation error in a generic application setting. For well-behaved systems with smooth basin structures, it may be possible to formalize the convergence characteristics of some resilience measures in terms of the sampling size of initial conditions. In Appendix \ref{app: synthetic}, we approach this endeavour in a simple system for which the true values of the resilience measures can be analytically derived. However, the convergence behaviour may be arbitrarily convoluted for other systems, such that an intelligent approach to grid tessellation and other numerical configurations becomes crucial. The low computational cost of our presented method enables such fine-tuning and trial-and-error procedures.

We anticipate that the framework presented here will be a useful basis for future research on existing and novel resilience concepts in applied and abstract contexts. Given that systems analysis with respect to resilience is a cornerstone of several interconnected research domains, this work constitutes an important advance towards a homogenization of terminology and methodology.

\newpage%%%%%%%%%%%%%%%%%%%%%%%%%%%%%

\begin{acknowledgments}
    This is ClimTip contribution \#98; the ClimTip project has received funding from the European Union's Horizon Europe research and innovation programme under grant agreement No. 101137601: Funded by the European Union. Views and opinions expressed are however those of the author(s) only and do not necessarily reflect those of the European Union or the European Climate, Infrastructure and Environment Executive Agency (CINEA). Neither the European Union nor the granting authority can be held responsible for them.
\end{acknowledgments}

\section*{Author Declarations}
\subsection*{Conflicts of Interest}
The authors have no conflicts to disclose.

\subsection*{Author Contributions}
AM: Conceptualization, Formal Analysis, Methodology, Software, Visualization, Writing - original draft and reviewing. CK: Conceptualization, Formal Analysis, Funding Acquisition, Methodology, Project Administration, Supervision, Writing - original draft and reviewing. GD: Conceptualization, Formal Analysis, Methodology, Software, Supervision, Visualization, Writing - original draft and reviewing.

\section*{Data Availability}
Visit the GitHub repository \href{https://github.com/andreasmorr/ResilienceComparison}{ResilienceComparison} to access the code generating all figures in this manuscript. The employed methodology can be used as part of the Julia package DynamicalSystems.jl. 
The implementation is covered by two names in the software: \verb|StabilityMeasuresAccumulator|
and \verb|stability_measures_along_continuation|.

\appendix
\mathtoolsset{showonlyrefs=false}
\section{Convergence properties}\label{app: synthetic}

Both the identification of system attractors as well as the estimation of the resilience measures associated with them are based on the sampling of initial conditions. In the first step, the trajectories corresponding to each initial condition are analysed for their convergence in the sense of recurrences in a tessellated grid, yielding the system attractors. For this, it should always be ensured that enough random initial conditions are considered in order for all relevant attractors to be found. In systems with smooth basin structures and only a handful of attractors, this is usually associated with permissible computational costs. When considering the evolution of attractors along a parameter change, the knowledge about the location of attractors in the previous parameter setting is used to make an educated guess about the sensible initial conditions for the next parameter value.

\begin{figure}[t]
	\centering
	\includegraphics[width=\columnwidth]{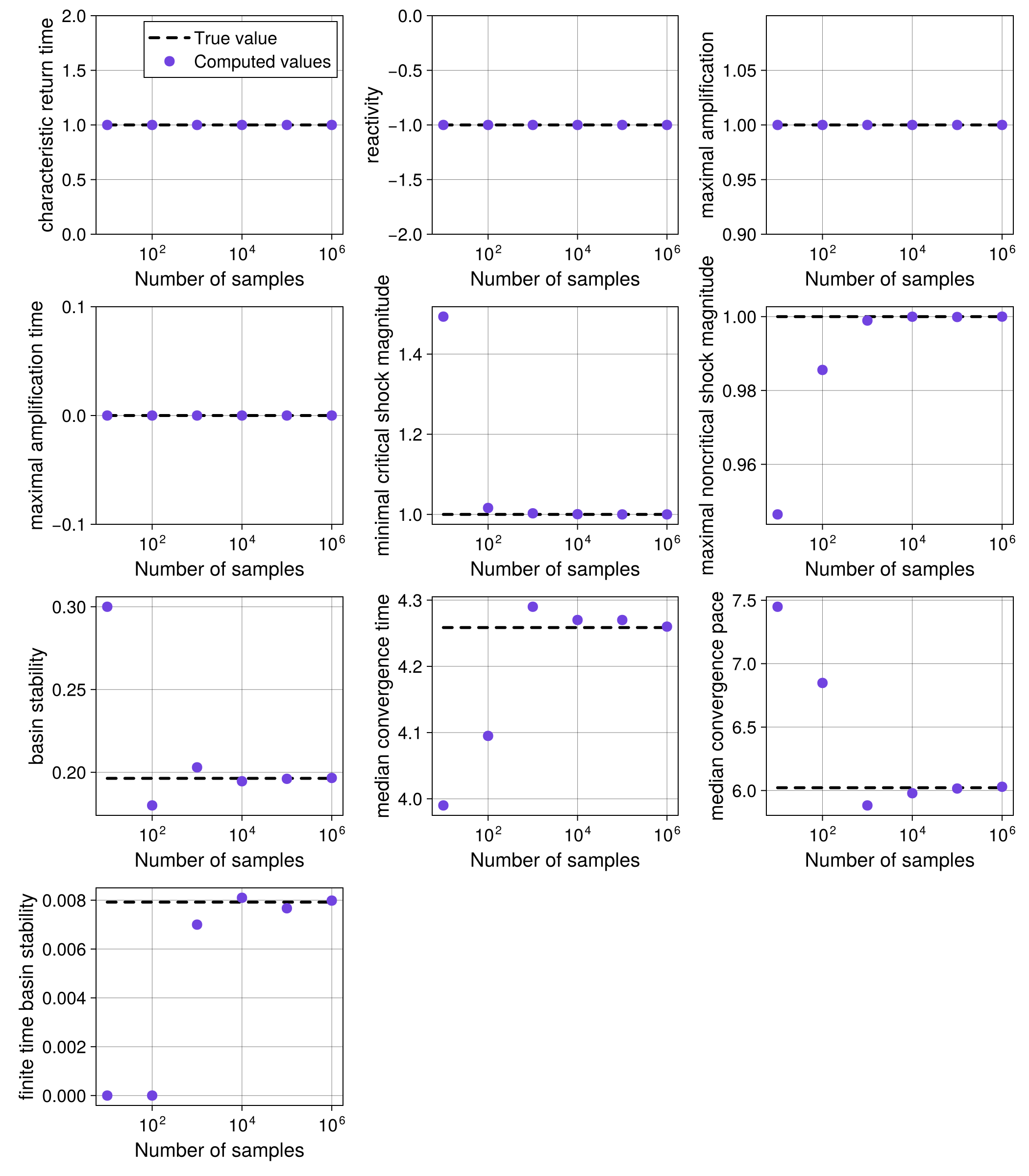}
	\caption{Estimated resilience measures of the point attractor in the system defined by Eq.~\ref{eq: simple}. They are plotted as a function of the employed initial condition sampling size in the respective estimation.}
    \label{fig: simple}
\end{figure}

The sample size for the second part of the procedure, i.e., estimating the resilience measures, directly determines the estimation error of some of the results. For the local resilience measures, this is not the case since they are merely a function of the estimated Jacobian matrix at the found attractors. For this, the functionality of ForwardDiff.jl is employed, which generally delivers satisfying accuracy. The nonlocal measures, however, directly depend on information gathered in the entire state space. For example, a satisfying estimation of the distance of the attractor to the nearest point not belonging to its own basin of attraction (minimal critical shock) necessitates a homogeneous and thorough exploration of the state space. It can therefore be expected that the estimated value of such measures approaches the true underlying values only for large initial condition samples. This is especially true for systems with convoluted or even fractal basin structures. For measures that rely on basin volume ((finite-time) basin stability), the random sampling of initial conditions essentially amounts to sampling a categorical distribution. The associated convergence of the sample fractions to the underlying probabilities is of order $\mathcal{O}(1/\sqrt{N})$, where $N$ is the sample size. \cite{Menck2013BasinStability} We generally approach the question of how the estimation of the resilience measures converges as a function of sample size in a simple system. This is described by the evolution equation
\begin{eqnarray}
    \frac{\dint x}{\dint t}&=\sigma(x^2+y^2-1)ax,\nonumber\\
    \frac{\dint y}{\dint t}&=\sigma(x^2+y^2-1)ay,\label{eq: simple}
\end{eqnarray}
where $\sigma$ is the sign function and $a>0$ determines the strength of attraction. Effectively, the system has a point attractor at $\textbf{0}$ whose basin of attraction is the unit circle. For this system, we can algebraically compute all resilience measures associated with the attractor. We then employ our framework of resilience estimation several times, each time with a different sample size of initial conditions. The results can be seen in Fig.~\ref{fig: simple}. For those measures that do depend on a homogeneous sampling of the state space, a convergence to the true underlying value can be observed.

\section{Code example}\label{app: code}

The following code example generates the data for Fig.~\ref{fig: predprey}.
It finds the attractors of the system, continues them along a parameter using global continuation, and estimates all resilience measures in parallel to the continuation.

% displayed code
\lstset{
  literate={ε}{$\varepsilon$}1
           {Δ}{$\Delta$}1
}
\begin{lstlisting}[language=Python]
# requires Attractors v1.29 or later
using DynamicalSystems

# define dynamical system
function predator_prey(u, p, t)
    A, B, C, D, E = Tuple(p)
    x, y = Tuple(u)
    s = x^2/(A*x^2 + B*x + 1)
    dx = x*(1 - x)*(x - E) - s*y
    dy = y*(C*s - D)
    return SVector(dx, dy)
end

p = [2.05, -2.6, 0.4, 1.0, 0.4]
u0 = [0.5, 0.02]

ds = CoupledODEs(predator_prey, u0, p, 
    diffeq = (alg = Vern9(), abstol = 1e-8, reltol = 1e-8)
)

# first pass: find attractors
# using recurrences
grid = (
    range(0, 1; length = 201), # x
    range(0, 0.05; length = 201), # y
)
mapper = AttractorsViaRecurrences(ds, grid;
    consecutive_recurrences = 1e5, 
    attractor_locate_steps = 1e4,
    consecutive_lost_steps = 1e6,
    Δt = 0.1, stop_at_Δt = true
)
ascm = AttractorSeedContinueMatch(mapper)

pidx = 5
prange = 0.35:0.003:0.45 # E range
pcurve = [[pidx => p] for p in prange]

sampler, = statespace_sampler(grid)

_, attractors_cont = global_continuation(
	ascm, prange, pidx, sampler; samples_per_parameter = 100
)

# second pass: estimate resilience
# with many more initial conditions
# and epsilon-based termination
result = stability_measures_along_continuation(
    ds, attractors_cont, pcurve, sampler; 
    ε = 0.001,
    finite_time = 100.0,
    samples_per_parameter = 1e5,
    proximity_mapper_options = (
        Δt = 0.1, stop_at_Δt = true
    )
)

# `result` contains all resilience measures
# for all parameter values
\end{lstlisting}

\bibliography{BibAll}

\end{document}